%% file: paper.tex
\documentclass[a4paper,11pt,twoside]{article}
\usepackage[left=3cm, right=3cm, top=3.5cm, bottom=3.5cm]{geometry}
\usepackage{type1cm}
\usepackage[utf8]{inputenc}
\usepackage[T1]{fontenc}
\usepackage{aecompl}
\usepackage{ae}
\usepackage{textcomp}
\usepackage[english]{babel}
\usepackage{graphicx}
\usepackage{pdfsync}
\usepackage{tikz}
\usetikzlibrary{arrows}
\usepackage{pgffor}
\usepackage{dsfont}
\usepackage{capt-of} 

\usepackage[leftcaption]{sidecap}
\usepackage{hyperref}
\usepackage{hypbmsec}
\usepackage{amsmath,amsthm,amsfonts}
\usepackage{bbm}
\usepackage{enumerate}

\theoremstyle{plain}
\newtheorem{thm}{Theorem}[section]
\newtheorem{prop}[thm]{Proposition}
\newtheorem{cor}[thm]{Corollary}
\newtheorem{lem}[thm]{Lemma}
\newtheorem{rem}[thm]{Remark}
\newtheorem{defn}[thm]{Definition}

\newcommand{\E}{\mathbb{E}}
\renewcommand{\P}{\mathbb{P}}
\newcommand{\Z}{\mathbb{Z}}
\newcommand{\N}{\mathbb{N}}
\newcommand{\R}{\mathbb{R}}
\newcommand{\C}{\mathcal{C}_2}
\renewcommand{\O}{\mathcal{O}}
\renewcommand{\d}{\mathrm{d}}

\newcommand{\ball}[1]{\mathcal{B}_{#1}}

\newcommand{\abs}[1]{\lvert#1\rvert}

\newcommand{\indicator}[1]{\mathds{1}_{#1}}

\usepackage{fancyhdr}
\hypersetup{%
  pdftitle = {Internal aggregation models on comb lattices},
  pdfauthor = {Wilfried Huss, Ecaterina Sava},
  pdfkeywords = {divisible sandpile, internal diffusion limited aggregation, rotor-router model}
}

\begin{document}
\title{Internal Aggregation Models on Comb Lattices}
\author{Wilfried Huss\footnote{Vienna University of Technology, Austria.}, 
Ecaterina Sava\footnote{Graz University of Technology, Austria.}}
\maketitle

\begin{abstract}
The two-dimensional comb lattice $\C$ is a natural spanning tree of the Euclidean lattice 
$\mathbb{Z}^2$. We study three related cluster growth models on $\C$: 
\emph{internal diffusion limited aggregation} (IDLA), 
in which random walkers move on the vertices of $\C$ until reaching an unoccupied  
site where they stop; \emph{rotor-router aggregation} in which particles perform deterministic walks, 
and stop when reaching a site previously unoccupied; and the \emph{divisible sandpile model} where at 
each vertex there is a pile of sand, for which, at each step, the mass exceeding $1$ is distributed 
equally among the neighbours. We describe the shape of the divisible sandpile cluster on $\C$, 
which is then used to give inner bounds for IDLA and rotor-router aggregation.
\end{abstract}



\section{Introduction}

Let $G$ be an infinite, locally finite and connected graph with a chosen
origin $o\in G$.  \emph{Internal diffusion limited aggregation (IDLA)} is a random walk-based
growth model, which was introduced by {\sc Diaconis and Fulton} \cite{diaconis_fulton_1991}.
In IDLA $n$ particles start at the origin of $G$, and each particle performs a simple
random walk until it reaches a vertex which was not previously occupied. There the particle
stops, and from now on occupies this vertex, and a new particle starts its journey at the origin.
The resulting random set of $n$ occupied sites in $G$ is called the \emph{IDLA cluster}, and will
be denoted by $A_n$.

IDLA has received increased attention in the last years. In 1992, \textsc{Lawler, Bramson and
Griffeath} \cite{lawler_bramson_griffeath} showed that for simple random walk
on $\Z^d$, with $d\geq 2$, the limiting shape of IDLA, when properly rescaled,
is almost surely an Euclidean ball of radius $1$. In 1995, \textsc{Lawler} 
\cite{lawler_1995} refined this result by giving estimates on the fluctuations. Recently
several improvements have been obtained. \textsc{Asselah and Gaudilli\`{e}re}
\cite{asselah_gaudilliere, asselah_gaudilliere_2} proved an upper bound
of order $\log(n)$ for the inner fluctuation $\delta_I$ and of order $\log^2(n)$
for the outer fluctuation $\delta_O$ in all dimensions $d\geq 2$. In \cite{asselah_gaudilliere_3}
they improve the upper bound on the inner fluctuation to $\sqrt{\log(\text{radius})}$,
for $d\geq 3$. Independently, and by different methods, \textsc{Jerison, Levine and Sheffield}
\cite{jerison_levine_sheffield,jerison_levine_sheffield2}, proved also that both $\delta_I$ and $\delta_O$ are of
order $\log(n)$ for IDLA on $\Z^2$ and of order $\sqrt{\log(\text{radius})}$ for $d\geq 3$ . 

\emph{Rotor-router aggregation} is a deterministic version of IDLA, where particles
perform \emph{rotor-router walks}, which are deterministic analogues to random walks. 
They have been first introduced into the physics literature under the name \emph{Eulerian walks} by 
\textsc{Priezzhev, D.Dhar et al}\cite{PhysRevLett.77.5079}.  At each vertex of the graph $G$,
we have an arrow (rotor) pointing to one of the neighbours of the vertex. A particle
performing a \emph{rotor-router walk} first changes the rotor at its current position to point
to the next neighbour, in a fixed order chosen at the beginning, and then moves to the neighbour
the rotor is now pointing at. 
In rotor-router aggregation each particle performs a rotor-router walk until it reaches an
unoccupied site, where it stops. Then a new particle starts at the origin, without resetting
the configuration of rotors. The resulting deterministic set $R_n$ of $n$ occupied sites is called
the \emph{rotor-router cluster}.

Rotor-router aggregation on the Euclidean lattice $\Z^d$ has been studied 
by \textsc{Levine and Peres} \cite{peres_levine_strong_spherical}, who showed that the cluster $R_n$ is
a ball in the Euclidean distance. On the homogeneous tree \textsc{Landau and Levine}
\cite{landau_levine_2009} proved that, provided the start configuration of rotors is acyclic,
the rotor-router cluster forms a perfect ball with respect to the
graph metric, whenever it has the right amount of particles. \textsc{Kager and Levine}
\cite{kager_levine_rotor_aggregation} studied the shape of the rotor-router cluster on a
modified two dimensional lattice, which they call the \emph{layered square lattice}.
In each of the known examples the limiting shape of rotor-router aggregation is the same
as the one for IDLA, but with much smaller fluctuations compared to IDLA.

In order to prove inner bounds for the above models, we use a third growth model,
the so-called \emph{divisible sandpile}, which has been introduced by \textsc{Levine and Peres} 
\cite{peres_levine_strong_spherical} as a tool for studying internal growth models
on $\Z^d$. In the divisible sandpile model each vertex can have an arbitrary amount
of mass. If a vertex has mass at least $1$, it is called \emph{unstable} and it can \emph{topple} 
by distributing the mass exceeding $1$ equally among its neighbours.  
At each timestep a vertex is chosen and toppled if it is unstable. Provided every vertex
is chosen infinitely often, the masses converge to a limiting distribution $\leq 1$.
The set of vertices with limit mass equal to $1$ is called the \emph{divisible sandpile cluster}.
If we start with a mass of $n$ concentrated at the origin, the corresponding sandpile cluster
will be denoted by $S_n$.

All three growth models have very similar behaviour. This was first noticed
by \textsc{Levine and Peres} \cite{peres_levine_strong_spherical,peres_levine_scaling_limits}
for the case when the state space is an Euclidean lattice. Computer
simulations suggest that the connection between the three growth models
holds in wide generality, but only partial results are available for
other state spaces. All these three models have the so-called
\emph{abelian property}, which makes them amenable to rigorous analysis.
In the case of IDLA and the rotor-router model this means that, if we let
several particles run at the same time, instead of one after another,
it is irrelevant for the end result in which order the particles make
their moves. In the case of the divisible sandpile model, it means that the
limiting distribution is independent of the order in which vertices topple. 

The aim of this paper is to study the three aggregation models described above 
on the \emph{comb lattice} $\C$, which is the spanning tree of the two-dimensional
Euclidean lattice $\Z^2$,  obtained by removing all horizontal edges of $\Z^2$ except the ones on
the $x$-axis.

\begin{minipage}{0.45\linewidth}
The graph $\C$ can also be constructed from a two-sided infinite path $\Z$
(the "\emph{backbone}" of the comb), by attaching copies of $\Z$ (the "\emph{teeth}") at
every vertex of the backbone. \\

We use the standard embedding of the comb into  $\Z^2$,
and use Cartesian coordinates $z = (x,y)\in \Z^2$ to denote vertices of $\C$. The vertex $o = (0,0)$
will be the \emph{root} or the \emph{origin}; see Figure \ref{fig:comb}. For functions $g$ 
on the vertex set of $\C$ we will often write $g(x,y)$ instead of $g(z)$, when $z = (x,y)$.
\end{minipage}
\hfill
\begin{minipage}{0.55\linewidth}
\centering
\input{comb}
\captionof{figure}{The comb $\C$}
\label{fig:comb}
\end{minipage}

While $\C$ is a very simple graph, it has some remarkable properties. For example,
the so-called \emph{Einstein relation} between the spectral-, walk- and fractal-dimension is
violated on the comb, see \textsc{Bertacchi} \cite{bertacchi_comb}. \textsc{Peres and Krishnapur}
\cite{peres_krishnapur_collide} showed that on $\C$ and other similar recurrent graphs two independent 
simple random walks meet only finitely often. Random walks on $\C$ have been 
studied by various authors, the first being \textsc{Havlin and Weiss} \cite{havlin_weiss_comb} 
and \textsc{Gerl} \cite{gerl86_zd_spanning_trees}.

The paper is organized as follows. In Section \ref{sec:preliminaries} we 
introduce some notations and basic facts which will be used through the rest of 
the work. Section \ref{sec:sandpile_model} is dedicated to the study of the divisible sandpile
on the comb $\C$. We show in Theorem \ref{thm:sandpile_cluster}
that the sandpile cluster $S_n$ on $\C$ has up to constant fluctuations the shape
\begin{align}\label{eq:limit_shape}
\ball{n} = \left\lbrace (x,y)\in\mathcal{C}_2:\: \frac{\abs{x}}{k} + 
\left(\frac{\abs{y}}{l}\right)^{1/2}\leq n^{1/3}\right\rbrace
\end{align}
where
\begin{equation*}
k = \left(\frac{3}{2}\right)^{2/3},
\quad l = \frac{1}{2}\left(\frac{3}{2}\right)^{1/3}.
\end{equation*}
Section \ref{sec:idla}
deals with IDLA on $\C$. Using the results obtained for the 
sandpile model, we prove an inner bound for IDLA, which is
of the type \eqref{eq:limit_shape}. Finally, in Section \ref{sec:rotor_router},
we give an inner estimate for the rotor-router model on $\C$ which is weaker than the
result obtained for IDLA. For a fixed initial configuration of rotors the exact shape of
the rotor-router cluster on the comb has been obtained by the authors in \cite{huss_sava_rotor}
using a purely combinatorial approach.

\section{Preliminaries}
\label{sec:preliminaries}

Let $(G,E(G))$ be an infinite, undirected and connected graph,
with vertex set $G$, equipped with a symmetric \emph{adjacency relation}
$\sim$, which defines the set of edges $E(G)$ (as a subset of $G\times G$). We 
write $(x,y)$ for the edge between the pair of neighbours $x,y$. In order to simplify 
the notation, instead of writing $(G,E(G))$ for a graph, we shall write only
$G$, and it will be clear from the context whether we are considering edges or vertices.
Let $o\in G$ be some fixed reference vertex called the \emph{origin}. 
For $x,y\in G$, let $\d(x,y)$ be the length of the shortest path from $x$ to $y$. 
Also, write $d(x)$ for the \emph{degree} of $x$, i.e, the number of neighbours of $x$.
For a subset $A\subset G$ we denote by
\begin{equation*}
\partial A = \bigl\{x\in G\setminus A:\: \exists y\in A \text{ with } x\sim y \bigr\}
\quad\text{and}\quad
\partial_I A = \bigl\{x\in A:\: \exists y\not\in A \text{ with } x\sim y \bigr\}
\end{equation*}
the \emph{(outer) boundary} respectively the \emph{inner boundary} of $A$.

Let $P=\big(p(x,y)\big)_{x,y\in G}$ be the one-step transition probabilities of the
simple random walk on $G$, i.e., $p(x,y)=1/d(x)$ if $y\sim x$ and $0$ otherwise.
We write $X_t$ for the position of the random walker at the discrete timestep $t$.
Probabilities will be written as $\P$, in particular $\P_x$ denotes the probability
of a random walk which starts at $x\in G$. Similarly $\E$ and $\E_x$ will denote
expectations using the same convention. For $y,z\in G$ the \emph{Green function} is defined as
\begin{equation*}
 G(y,z)=\mathbb{E}_y\Big[\sum_{t=0}^\infty\mathbf{1}_{\{X_t=z\}}\Big],
\end{equation*}
and represents the expected number of visits to $z$ of the
random walk $X_t$ started at $y$. For a subset $A\subset G$, write $G_A$ for
the Green function of the random walk stopped upon leaving the set $A$. That is,
if $\tau=\min\{t\geq 0: X_t\notin A\}$, then 
\begin{equation*}
G_A(x,y)= \mathbb{E}_x\Big[\sum_{t=0}^{\tau-1}\indicator{\{X_t=y\}}\Big].
\end{equation*}
For a function $f:G\to \mathbb{R}$, its \emph{Laplace operator} $\Delta f$ is defined as
\begin{equation*}
\Delta f(x) = \frac{1}{d(x)}\sum_{y\sim x}\big( f(y)-f(x)\big).
\end{equation*}
A function $f:G\to\mathbb{R}$ is called \emph{superharmonic} on a set $A\subset G$ if $\Delta f\leq 0$,
and \emph{harmonic} if $\Delta f=0$, for all $x\in A$.  For a function $g: G\rightarrow \R$, define
its \emph{least superharmonic majorant} as
\begin{equation*}
s(x) = \inf \big\lbrace f(x):\: f \text{ superharmonic } ,\, f \geq g \big\rbrace.
\end{equation*}
Remark that the function $s$ is itself superharmonic on $G$. The following is widely known.
\begin{lem}[Minimum principle.] If $f$ is a superharmonic function on $G$ 
and there exists $x\in G$ such that $f(x)=\min_{G}f$, then $f$ is constant.
\end{lem}
\section{Divisible Sandpile}\label{sec:sandpile_model}

Let $\C$ be the comb as in Figure \ref{fig:comb}, and let $\mu_0$ be an 
{\em initial mass distribution} on $\C$, i.e., a function 
$\mu_{0}:\C\to\mathbb{R}_{+}$ with finite support.
The \emph{divisible sandpile} is a sequence $(\mu_k)_{k\geq 0}$ of  
mass distributions, which are created according to the following rule. 
At each time step $k$, choose a vertex $x\in\C$. If $\mu_k(x)\geq 1$, 
the pile of sand at $x$ is unstable and topples, which means that 
$x$ keeps mass $1$ for itself and the remaining mass $\mu_k(x)-1$ 
is distributed equally among the neighbours $y$ of $x$, that is, according to 
the transition probabilities $p(x,y)$ of the simple random walk on $\C$.
Given a mass distribution $\mu_k$ at time $k$ and a vertex $x\in\C$, the 
\emph{toppling operator} can be defined as
\begin{equation*}
T_x \mu_k(y) = \mu_k(y) + \alpha_k(x) d(y) \Delta\delta'_x(y),  \text{ for } y\in \C,
\end{equation*}
where $\delta'_x(y)=\frac{\delta_x(y)}{d(y)}$ and $\alpha_k(y) = \max \lbrace \mu_k(y)-1, 0 \rbrace$.
Let $(x_k)_{k\geq 0}$ be a sequence of vertices in $\C$ called 
the \emph{toppling sequence}, which contains each vertex of $\C$ infinitely often. 
Then  the mass distribution of the sandpile after $k$ steps is defined as
\begin{equation*}
\mu_{k+1} = T_{x_k}\mu_k = T_{x_k}\cdots T_{x_0}\mu_0.
\end{equation*}
Hence, $\mu_{k+1}(y)$ is  the amount of mass present at $y$ after toppling the sites 
$x_0,\ldots,x_k$ in succession.
One of the tools that will be used throughout this work in various incarnations 
is the so-called \emph{odometer function}, which was introduced by 
\textsc{Levine and Peres}
\cite{peres_levine_strong_spherical}.
\begin{defn} 
The \emph{odometer function} $v_k$ is defined as
\begin{equation*}
v_k(y) =  \sum_{j\leq k:\:x_j=y} \mu_j(y)-\mu_{j+1}(y) = \sum_{j\leq k:\:x_j=y}\alpha_j(y),\ y\in\C,
\end{equation*}
and represents the total mass emitted from $y$ during the first $k$ topplings.
\end{defn}
For simple random walks on $\C$ it is easier to work with the \emph{normalized odometer function} 
$u_k(x)=\frac{v_k(x)}{d(x)}$. Lemma $3.1$ of {\sc Levine and Peres} \cite{peres_levine_strong_spherical}
can be easily adapted to our case, in order to show that, as $k$ goes to infinity,
$\mu_k$ and $u_k$ converge to limit functions $\mu$ and $u$ respectively. Define 
\begin{equation*}
S=\lbrace x\in \C:\: \mu(x) = 1 \rbrace.
\end{equation*}
The set $S$ is called the \emph{sandpile cluster} with initial mass distribution $\mu_0$.
The limit functions $\mu$ and $u$ satisfy
\begin{equation}\label{eq:mu_odometer}
 \mu(x)=\mu_0(x)+d(x)\Delta u(x),\text{ for all }x\in S,
\end{equation}
and 
\begin{equation}\label{eq:mu_leq_1}
 \mu(x)\leq 1, \text{ for all }x\in \C.
\end{equation}
From \eqref{eq:mu_odometer} and \eqref{eq:mu_leq_1} it follows that
\begin{equation}
\label{eq:sandpile_odometer_laplace}
\Delta u(x) = \frac{1}{d(x)}\bigl(1 - \mu_0(x)\bigr) \quad \text{for all } x\in S,
\end{equation}
and $u(x) = 0$, if $x\not\in S$.
The following result provides a method for solving this free boundary problem.
For a proof, see once again {\sc Levine and Peres} \cite[Lemma $3.2$]{peres_levine_strong_spherical}.
\begin{lem}
\label{lem:calc_odometer}
Consider a function $\gamma:\C\rightarrow \R$ with
\begin{equation}
\label{eq:odometer_potential}
\Delta \gamma(x)  =  \dfrac{1}{d(x)}\bigl(1-\mu_0(x)\bigr), \quad \text{for all } x\in \C.
\end{equation}
Then the normalized odometer function $u$ of the sandpile satisfies $u = \gamma + s$, 
where $s$ is the least superharmonic majorant of $-\gamma$.
\end{lem}
Lemma \ref{lem:calc_odometer} gives a representation of the odometer function which is
independent of the toppling sequence. 
\begin{rem}[Abelian property]
\label{rem:abelian_property}
The limit $u$ of the normalized odometer function and the sandpile cluster $S$ are independent
of the toppling sequence  $(x_k)_{k\geq 0}$.
\end{rem}

\subsection{Divisible Sandpile on the Comb}

With the help of Lemma \ref{lem:calc_odometer}, we shall next describe 
the limit shape of the sandpile cluster on the two-dimensional
comb $\C$. Consider an initial mass distribution
$\mu_0$ concentrated at the origin $o$, that is $\mu_0=n\cdot\delta_o$, and denote by 
\begin{equation*}
S_n=\{z\in\C:\mu(z)=1\}
\end{equation*}
the sandpile cluster, and by $u_n$ the limit of the normalized odometer  function for this choice of initial distribution. 
We use another simple fact about $u_n$; for a proof see \cite[Lemma 3.4]{peres_levine_strong_spherical}.
\begin{lem}
\label{lem:u_strictly_decreasing}
If $x\in S_n\setminus\{o\}$ and $y\sim x$ with $\d(o,y) < \d(o,x)$, then $u_n(y) \geq u_n(x) + 1$.
\end{lem}
By \eqref{eq:sandpile_odometer_laplace}, the normalized odometer function satisfies
\begin{equation}\label{eq:laplace_odometer_function}
\Delta u_n(z)  =  \dfrac{1}{d(z)}\bigl(1-n\cdot \delta_o(z)\bigr), \text{ for }z\in S_n.
\end{equation}
The odometer function $u_n$ can be reduced to odometer 
functions of suitable divisible sandpiles on $\Z$, which are easy to compute.
Let $\tilde{u}_n$ be the normalized odometer function of the divisible sandpile on 
$\Z$, with initial mass distribution $\tilde{\mu}_0$ concentrated at $0$, that is, 
$\tilde{\mu}_0 = n \cdot\delta_0$. By Remark \ref{rem:abelian_property} it is 
clear that the sandpile cluster $\tilde{S}_n$ on $\Z$ in this case is a 
symmetric interval around the origin $0$.
In order to compute $\tilde{u}_n$, by Lemma \ref{lem:calc_odometer},
we need to construct a function $\tilde{\gamma}_n:\Z\to\R$ with Laplacian given in \eqref{eq:odometer_potential}.
It is easy to check that $\tilde{\gamma}_n$ defined by
\begin{equation}
\label{eq:gamma_z}
\tilde{\gamma}_n(y) = \frac{1}{2}\left(|y| - \frac{n}{2}\right)^2,
\end{equation}
satisfies the required property. Since $\tilde{\gamma}_n$ is nonnegative, the constant function $0$ is a superharmonic
majorant of $-\tilde{\gamma}_n$. Hence, by Lemma \ref{lem:calc_odometer}, we have $\tilde{u}_n \leq \tilde{\gamma}_n$.
Now, consider $\gamma_n:\C\to \R$ with
\begin{equation}\label{eq:odometer_C2}
\gamma_n(x,y) = \tilde{\gamma}_{n_x}(y),\text{ for }(x,y)\in\C,
\end{equation}
where $n_x\in\R$ for all $x\in \Z$. The quantities $n_x$ can be interpreted as the total amount of 
mass that ends up in the copy of $\Z$ that is attached to the vertex $(x,0)$. 
Then $\gamma_n$ satisfies \eqref{eq:odometer_potential} if and only if
\begin{align}
\label{eq:odometer_comb_inc_dimension}
n_x = n\cdot\indicator{\{x=0\}} + \tilde{\gamma}_{n_{x-1}}(0) - 2\tilde{\gamma}_{n_x}(0)
      + \tilde{\gamma}_{n_{x+1}}(0)
\end{align}
holds for all $x\in\Z$. From \eqref{eq:gamma_z} and \eqref{eq:odometer_comb_inc_dimension}, and using
the fact that $n_x = n_{-x}$ by symmetry, we get the following recursion for the numbers $n_x$
\begin{align}
\label{eq:odometer_comb_recursion_1}
n_0 &= n + \frac{1}{4} n_1^2 - \frac{1}{4} n_0^2, \\
\label{eq:odometer_comb_recursion_2}
n_x &= \frac{1}{8} n_{x-1}^2 - \frac{1}{4} n_x^2 + \frac{1}{8} n_{x+1}^2, \text{ for } x > 0.
\end{align}

\begin{figure}[t]
\centering
\includegraphics[height=7cm]{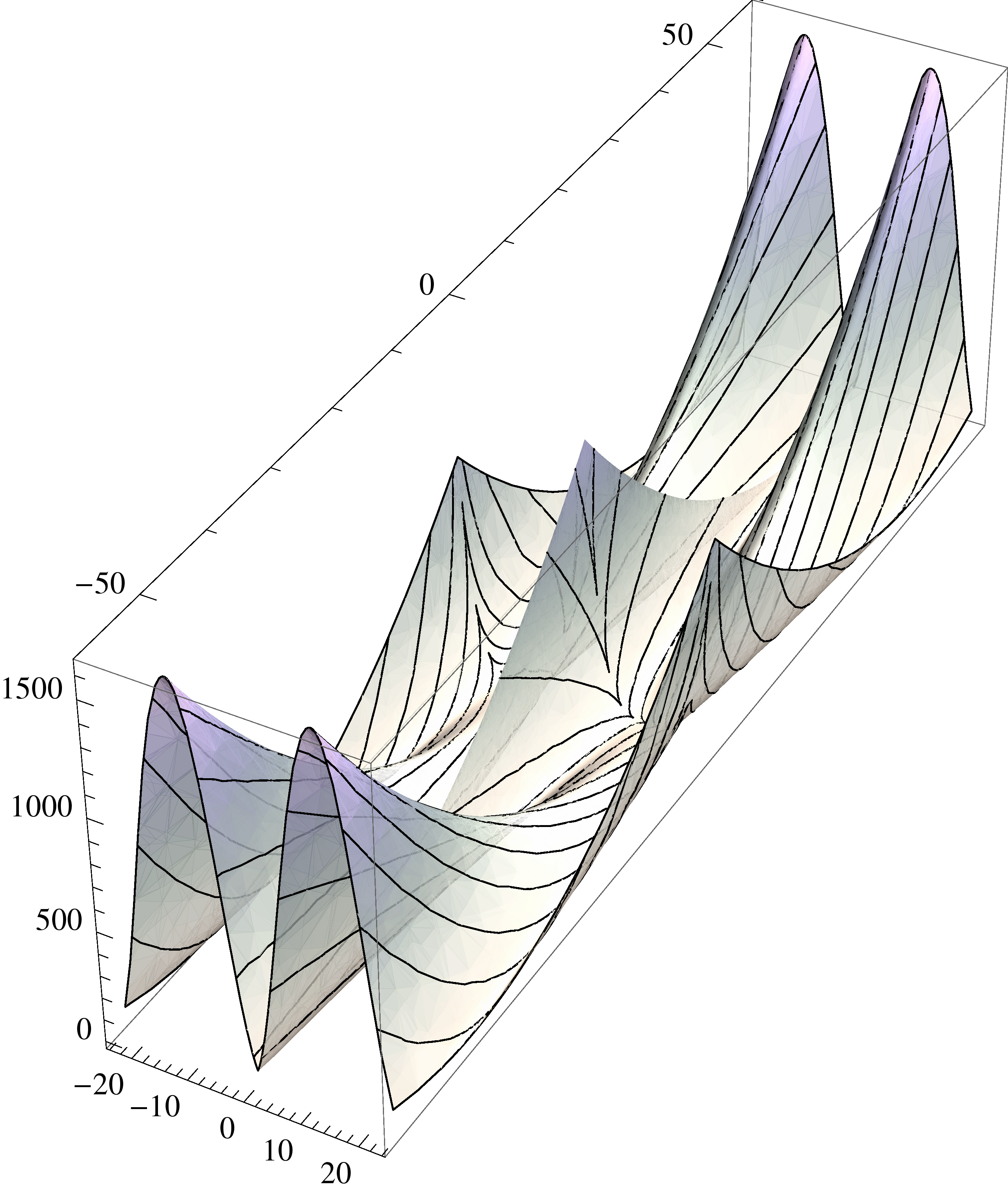} \hspace{3cm}
\includegraphics[height=7cm]{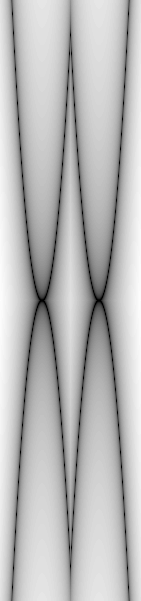}
\caption{\label{fig:gamman} Two plots of $\gamma_n$ for $n=1000$. The
graphic on the left is superimposed with contour lines representing the sets
$\ball{n}$ for various values of $n$. In the density plot on the right, dark areas
represent small values. By construction, the finite area which is surrounded by the
local minima of $\gamma_n$ coincides with the region $S_n$ covered by the sandpile.}
\end{figure}

Equation \eqref{eq:odometer_comb_recursion_2} has strictly positive solutions as quadratic polynomials of the form
\begin{equation}
\label{eq:n_x}
n_x = \frac{2}{3} x^2 - t\cdot x + \frac{9 t^2 + 4}{24},  \text{ with }t\in\R.
\end{equation}
By the initial condition \eqref{eq:odometer_comb_recursion_1}, the parameter $t$ satisfies the equation
\begin{equation*}
n = \frac{3}{16} t^3 + \frac{5}{12} t,
\end{equation*}
which has one real root given by
\begin{equation}
\label{eq:t}
t = T(n) - \frac{20}{27} T(n)^{-1},
\end{equation}
with
$T(n) = \left(\frac{8\sqrt{3}}{243}\sqrt{2187n^2 + 125} + \frac{8}{3}n \right)^{\frac{1}{3}}$.
By a series expansion around $n=\infty$, one obtains
\begin{equation}\label{eq:t_series}
t = 2 \left(\frac{2}{3}\right)^{1/3} n^{1/3} + \O(1).
\end{equation}
Therefore, the function $\gamma_n(x, y) = \tilde{\gamma}_{n_x}(y)$
with $n_x$ defined by \eqref{eq:n_x} and \eqref{eq:t} satisfies the conditions of 
Lemma \ref{lem:calc_odometer}, with $\mu_0=n\cdot\delta_o$. See Figure \ref{fig:gamman}
for a graphical representation of $\gamma_n$. We are now ready to prove the limit shape 
for the divisible sandpile on $\C$.

\begin{thm}
\label{thm:sandpile_cluster}
Let $S_n$ be the divisible sandpile cluster on $\C$, with $\mu_0 = n\cdot\delta_o$. Then there exists 
a constant $c\geq 0$ such that, for $n \geq n_0$:
\begin{align*}
\ball{n-c} \subset S_n \subset \ball{n+c},
\end{align*}
where
\begin{align}\label{eq:lim_cluster}
\ball{n} = \left\lbrace (x,y)\in\mathcal{C}_2:\: \frac{\abs{x}}{k} + 
\left(\frac{\abs{y}}{l}\right)^{1/2}\leq n^{1/3}\right\rbrace
\end{align}
and
\begin{equation*}
k = \left(\frac{3}{2}\right)^{2/3}, \quad \quad 
l = \frac{1}{2}\left(\frac{3}{2}\right)^{1/3}.
\end{equation*}

\end{thm}

\begin{proof} \emph{The upper bound $S_n\subset \ball{n+c}$:}
The mass distributions $n_x$ are nonnegative
for all $x$, therefore $\gamma_n$ is nonnegative, and this implies that
the constant function $0$ is a superharmonic majorant of $-\gamma_n$. 
Thus, by Lemma \ref{lem:calc_odometer}, $\gamma_n$ is an upper 
bound of the odometer function $u_n$. Moreover, Lemma \ref{lem:u_strictly_decreasing}
implies that $u_n$ decreases by a fixed amount on the sandpile cluster $S_n$ when we move away from the origin.
Therefore, in order to get an upper bound for $S_n$, it suffices to calculate the minima of $\gamma_n$ along each
infinite ray starting at $o=(0,0)$. By the symmetry of $\C$, it is sufficient to consider only the first quadrant.

Consider the rays which lies entirely on the positive $x$-axis. We have
$\gamma_n(x,0) = \frac{1}{8} n_x^2$. The minimum of this function 
is attained at $x_{\min}=\frac{3}{4}t$, with $t$ given in \eqref{eq:t}. 
Using the series expansion \eqref{eq:t_series} of $t$ we get
\begin{align}
\label{eq:x_min}
x_\text{min} = k n^{1/3} + \O(1), \quad \text{ with }k = \left(\frac{3}{2}\right)^{2/3},
\end{align}
which is also an upper bound of $S_n$ on the $x$-axis by Lemma
\ref{lem:u_strictly_decreasing}, since $\gamma_n(\lfloor x_{\min}\rfloor, 0)$ is bounded
by a constant which is independent of $n$, and smaller than $1/10$.

To calculate the extent of the sandpile cluster on the ``teeth'', we need to compute the
minima of $\gamma_n$ in the $y$-direction. On each ``tooth'' of the comb, $\gamma_n$ is a
quadratic polynomial which attains its minimum at $y_{\min}(x) = \frac{n_x}{2}$. Moreover,
$\gamma\big(x,\lfloor y_{\min}(x) \rfloor\big)\leq 1/2$. 
Using \eqref{eq:n_x} and a series expansion around infinity we get
\begin{equation*}
y_\text{min}(x) 
= l\left(n^{1/3} - \frac{x}{k}\right)^2 + \frac{2}{3} x -
 \frac{1}{2l} n^{1/3} - \frac{7l}{9k} x n^{-1/3} +\mathcal{O}(1),
\end{equation*}
where $l = \frac{1}{2}\left(\frac{3}{2}\right)^{1/3}$.
By the estimate in the $x$-direction we know that $(x,y)\in S_n$ only if 
$x\leq x_{\min}$.
Thus, using the expansion \eqref{eq:x_min} for $x_{\min}$ we obtain $(x,y) \in S_n$ if
$\abs{x} \leq k n^{1/3} +\O(1)$ and
$\abs{y} \leq l\left(n^{1/3} - \frac{x}{k}\right)^2+\O(1)$, for $n\geq n_0$.
This proves the upper bound $S_n \subset \ball{n+c}$.

\emph{The lower bound $\ball{n-c} \subset S_n$:} On each infinite ray the
minimum of $\gamma_n(z)$ is smaller than a constant $a>0$, independent of $n$.
Also, from the upper bound, we have $u_n(z) = 0$ for all $z\in\partial \ball{n+c}$.
Hence $u_n(z) - \gamma_n(z) \geq -a$ for all $z\in\partial \ball{n+c}$. Since the
function $u_n - \gamma_n$ is superharmonic, by the \emph{Minimum Principle}, it attains its
minimum on the boundary and the inequality
$u_n(z) - \gamma_n(z) \geq -a$
holds for all $z\in\ball{n+c}$. Thus $\gamma_n - a$ is also a lower bound of the odometer
function on $\ball{n+c}$, which gives the inner estimate
$\ball{n-c}\subset S_n$, for some constant $c$.
\end{proof}
The next corollary follows directly from the proof of the theorem.
\begin{cor}
\label{cor:sandpile_odometer}
Let $u_n$ be the normalized odometer function of the divisible sandpile on
$\C$, with initial mass distribution $\mu_0 = n\cdot\delta_o$, and
$\ball{n}\subset\C$ defined as in \eqref{eq:lim_cluster}.
Then there exists a constant $0<a<2$, such that, for all $n>n_0$ and all $z\in\C$
\begin{equation*}
\big(\gamma_n(z) - a\big)\indicator{\ball{n}}\leq u_n(z) \leq \gamma_n(z).
\end{equation*}
\end{cor}
\section[IDLA]{Internal Diffusion Limited Aggregation}
\label{sec:idla}

Let $\big( X^i_t \big)_{i\in\N}$ be a sequence of independent and identically
distributed simple random walks on the comb $\C$, with common starting point
$X^i_0 = o$. Then $X_t^i$ represents the position of the $i$-th particle
at time $t$. 
\begin{defn}
\label{def:idla}
\emph{Internal diffusion limited aggregation (IDLA)} is a stochastic
process of increasing subsets $\big( A_i \big)_{i\in\N}$ of $\C$, which
are defined recursively as $A_1 = \{ o \}$ and for $i\geq 2$
\begin{align*}
\P\big[A_i = \mathbf{A} \cup \lbrace x \rbrace |\: A_{i-1} =\mathbf{A}\big]
     = \P\big[X^{i}_{\sigma^{i}} = x\big],
\end{align*}
where $\sigma^{i} = \inf \lbrace t\geq 0:\: X^i_t \not\in A_{i-1} \rbrace$ is
the first exit time of the random walk $X^i_t$ from $A_{i-1}$.
\end{defn}

\begin{minipage}{0.55\linewidth}
The IDLA cluster is build up one site at a time. That is,
suppose that we already have the cluster $A_{i-1}$ 
after $i-1$ particles stopped, and we want to get $A_i$.
For this, the $i$-th particle $X^i_t$ starts at  $o$, and evolves as long it stays 
inside the IDLA-cluster $A_{i-1}$. When $X^i_t$
leaves $A_{i-1}$ for the first time, it stops, and the point outside of the cluster that is visited by $X^i_t$ 
is added to the new cluster $A_{i}$. The set $A_i$ is called the \emph{IDLA-cluster} of $i$ particles.

Figure \ref{fig:idla} shows IDLA clusters on $\C$ with $100, \ 500$ and $1000$ particles.\\

We will prove the following shape theorem for IDLA on  $\mathcal{C}_2$.
\end{minipage}
\hspace*{0.5cm}
\begin{minipage}{0.40\linewidth}
\centering
\input{idla_comb}

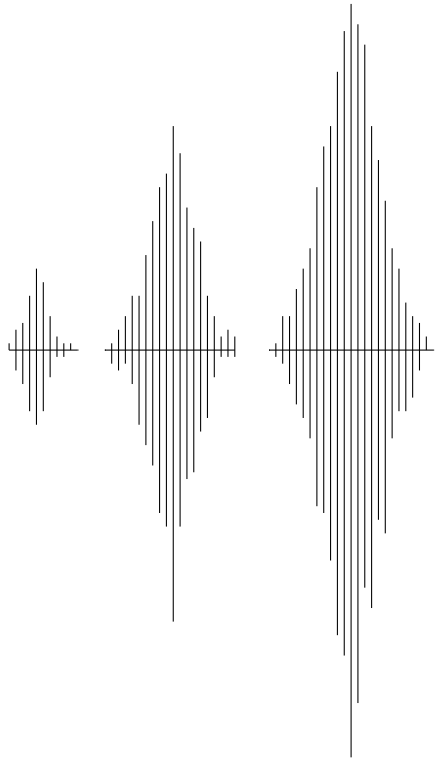
\captionof{figure}{IDLA cluster}
\label{fig:idla}
\end{minipage}
\begin{thm}\label{thm:idla_cluster}
Let $A_n$ be the IDLA cluster of $n$ particles on $\C$. 
Then, for all $\varepsilon >0$, we have with probability $1$ 
\begin{align}\label{eq:idla_inner_bound}
\ball{n(1-\varepsilon)} \subset A_n,\text{ for all sufficiently large } n,
\end{align}
where
\begin{align*}
\ball{n} = \left\lbrace (x,y)\in\mathcal{C}_2:\: \frac{\abs{x}}{k} + 
\left(\frac{\abs{y}}{l}\right)^{1/2}\leq n^{1/3}\right\rbrace
\end{align*}
and
\begin{equation*}
k = \left(\frac{3}{2}\right)^{2/3},
\quad l = \frac{1}{2}\left(\frac{3}{2}\right)^{1/3}. 
\end{equation*}
\end{thm}
The set $\ball{n}$ is the same as the limit shape of the divisible sandpile from Theorem \ref{thm:sandpile_cluster}.
The proof of Theorem \ref{thm:idla_cluster} uses ideas of \textsc{Lawler, Bramson and Griffeath} 
\cite{lawler_bramson_griffeath} and of \textsc{Levine and Peres} \cite{peres_levine_scaling_limits}. 
Following \cite{lawler_bramson_griffeath}, we introduce the stopping times
\begin{align*}
\tau^i_{n} = \min\{t\geq 0:X_t^i\notin \ball{n}\} \quad \text{ and }
 \tau^i_{z} = \min\{t\geq 0:X_t^i=z\},
\end{align*}
for $z\in\ball{n}$. Consider the probability that a fixed vertex $z\in\ball{n}$
does not belong to the IDLA cluster $A_n$, which can be written in terms of the stopping times
defined above as
\begin{equation*}
 \mathbb{P}[z\notin A_n]=\mathbb{P}\Big[\bigcap_{i\leq n}\sigma^i<\tau_z^i\Big].
\end{equation*}
Hence, by the Borel-Cantelli Lemma, convergence of the series
\begin{equation}\label{eq:borel_cantelli_idla}
 \sum_{n\geq n_0}\sum_{z\in \ball{n(1-\varepsilon)}}\mathbb{P}[z\notin A_n],
\end{equation}
is a sufficient condition for Theorem \ref{thm:idla_cluster}.
Fix now $n$ and $z\in \ball{n}$ and consider the random variables
\begin{equation*}
N=\sum_{i=1}^n\indicator{\{\tau^i_z<\sigma^i\}}, \quad
M=\sum_{i=1}^n\indicator{\{\tau^i_z<\tau^i_n\}} \quad \text{and} \quad
L=\sum_{i=1}^n\indicator{\{\sigma^i\leq\tau_z^i<\tau_n^i\}}
\end{equation*}
Then $N$ represents the number of particles that visit 
$z$ before leaving the cluster. The variable $M$ counts the number of particles
that visit $z$ before leaving $\ball{n}$, and $L$ is the number of particles that
visit $z$ after leaving the cluster $A_i$
but before leaving $\ball{n}$. Remark that if $L<M$, then $z\in\C$ belongs to  $A_n$. 
Moreover $N\geq M-L$. Therefore, in order to estimate $\mathbb{P}[z\notin A_n]$, 
we just need an upper bound for $\mathbb{P}[M=L]$. For any fixed number $a$
\begin{align}\label{eq:prb_M_L_a}
\begin{split}
\mathbb{P}[z\notin A_n] &= \mathbb{P}[N=0]\leq \mathbb{P}[M-L=0]\leq\mathbb{P}[M\leq a \text{ or }L\geq a]\\
		     &\leq  \mathbb{P}[M\leq a] +\mathbb{P}[L\geq a].
\end{split}
\end{align}
We shall show that for a suitable choice of $a$, the probabilities $\mathbb{P}[M\leq a]$ and 
$\mathbb{P}[L\geq a]$ are small enough, such that the series
\eqref{eq:borel_cantelli_idla} converges. The derivation of a suitable value of $a$ will be done differently,
not like in the case of Euclidean lattices, studied by \textsc{Lawler, Bramson and Griffeath} in
\cite{lawler_bramson_griffeath}, who used classical asymptotics for the Green function stopped on a ball.
Since in our case, the Green function $G_{\ball{n}}$ stopped on $\ball{n}$ is not directly available,
we use the odometer function of the divisible sandpile as a replacement,
as suggested by \textsc{Levine and Peres} in \cite{peres_levine_scaling_limits}.
For simplicity of notation, we shall write $G_n(y,z)$ instead of $G_{\ball{n}}(y,z)$.
The random variable $M$ is a sum of i.i.d. indicator variables, with
\begin{equation}\label{eq:expectation_M}
 \mathbb{E}[M]=n \mathbb{P}_{o}[\tau_z<\tau_n]=n\frac{G_n(o,z)}{G_n(z,z)}.
\end{equation}
Even though $L$ is a sum of dependent indicator variables, following \cite{lawler_bramson_griffeath},
$L$ can be bounded by a sum of independent indicators
as follows. Only those particles with $X^i_{\sigma^i}\in\ball{n}$ contribute
to $L$ and for each $y\in\ball{n}$ there is at most one index $i$ 
with $X^i_{\sigma^i}=y$. The corresponding post-$\tau_y$ random walks are
independent. In order to avoid dependencies in $L$, enlarge the 
index set to all of $\ball{n}$ and define
\begin{equation*}
 \tilde{L}=\sum_{y\in\ball{n}}\indicator{\{\tau_z<\tau_n\}}^y,
\end{equation*}
where the indicators $\indicator{}^y$ correspond to independent random walks
starting at $y$. Then $L\leq\tilde{L}$, and the
expectation of $\tilde{L}$ is given by
\begin{equation}\label{eq:expectation_L_tilde}
\mathbb{E}[\tilde{L}]=\sum_{y\in\ball{n}}\mathbb{P}_y[\tau_z<\tau_n]=\frac{1}{G_n(z,z)}\sum_{y\in \ball{n}}G_n(y,z). 
\end{equation}
Now \eqref{eq:prb_M_L_a} can be rewritten as
\begin{equation}\label{eq:prb_M_tilde_L_a}
\mathbb{P}[z\notin A_n] \leq \mathbb{P}[M\leq a] +\mathbb{P}[\tilde{L}\geq a]. 
\end{equation}
We shall relate the random variables $\tilde{L}$ and $M$ with the odometer function
of the divisible sandpile. For this, consider the function $f_n:\partial \ball{n}\cup\ball{n}\to \mathbb{R}$,
\begin{equation}\label{eq:fn_z}
 f_n(z)=\frac{G_n(z,z)}{d(z)}\mathbb{E}[M-\tilde{L}].
\end{equation}
Set $h_x(y) = \frac{G_n(x,y)}{d(y)}$, then for $x,y\in\ball{n}$ one has
\begin{equation}\label{eq:laplace_hy}
\Delta h_x(y) =
\begin{cases}
-\dfrac{1}{d(x)}, \quad & \text{ if } x=y \\
0, & \text{ if } x\neq y
\end{cases},
\end{equation}
By linearity of the Laplace operator, $f_n$ solves the following Dirichlet problem
\begin{equation*}
\begin{cases}
\Delta f_n(z) = \dfrac{1}{d(z)}\big(1-n\cdot\delta_o(z)\big), \quad & \text{ for } z\in\ball{n} \\
f_n(z) = 0, & \text{ for } z\in\partial\ball{n}.
\end{cases}
\end{equation*}
Recall that the odometer function $u_n$ of the divisible sandpile with initial distribution 
$\mu_0=n\cdot\delta_{o}$ solves the same Dirichlet problem \eqref{eq:laplace_odometer_function}
on the sandpile cluster $S_n$ (whose shape is given by $\ball{n}$ by Theorem \ref{thm:sandpile_cluster}).

The uniqueness of the solution of a Dirichlet problem gives that $f_n=u_n$ on the set $\ball{n}$, and $u_n$ is
approximated (up to an additive constant) by the function $\gamma_n$ defined in \eqref{eq:odometer_C2}. 
Since $u_n>0$, it follows that $f_n(z)>0$, for all $z\in\ball{n}$, which is equivalent to
$\mathbb{E}[M]>\mathbb{E}[\tilde{L}]$.

We will use the following large deviations estimate for sums of independent
indicators. For a proof, see \textsc{Alon and Spencer} \cite[Cor. A.1.14]{alon_spencer}.
\begin{lem}\label{lem:large_dev}
If $N$ is a sum of finitely many independent indicator random variables,
then for all $\lambda>0$,
\begin{equation*}
\mathbb{P}\big[|N-\mathbb{E}N|>\lambda \mathbb{E}N\big]<2 e^{-c_{\lambda}\mathbb{E}N},
\end{equation*}
where $c_{\lambda}$ is a constant depending only on $\lambda$.
\end{lem}
In order to find an upper bound for the right hand-side of \eqref{eq:prb_M_tilde_L_a}
we use the previous Lemma and choose $\lambda>0$ and $a$ such that
\begin{equation}\label{eq:ineq_a}
(1+\lambda)\mathbb{E}[\tilde{L}]\leq a \leq (1-\lambda)\mathbb{E}[M].
\end{equation}
Hence $\lambda$ has to satisfy the relation
\begin{equation}\label{eq:lambda}
 0<\lambda\leq\dfrac{\mathbb{E}[M-\tilde{L}]}{\mathbb{E}[M+\tilde{L}]} = \frac{f_n(z)}{g_n(z)},
\end{equation}
and $g_n$ defined as
\begin{equation}\label{eq:gn_z}
 g_n(z)=\frac{G_n(z,z)}{d(z)}\mathbb{E}[M+\tilde{L}].
\end{equation}
To obtain $\lambda$, we have to bound $f_n(z)/g_n(z)$  away from $0$.
For this, we first have to calculate the function $g_n$ which is,
like $f_n$, defined on $\partial \ball{n}\cup \ball{n}$.

\subsection{The function $g_n$}
By \eqref{eq:expectation_M} and \eqref{eq:expectation_L_tilde}, the function $g_n$ is the solution of
the Dirichlet Problem
\begin{equation}
\label{eq:g_laplace}
\begin{cases}
\Delta g_n(z) = \dfrac{1}{d(z)}\big(-1-n\cdot\delta_o(z)\big), \quad & \text{for } z\in\ball{n}, \\
g_n(z) = 0, & \text{for } z\in\partial\ball{n},
\end{cases}
\end{equation}
and can therefore be obtained by solving some linear recursions.
For simplicity, we first shift the set $\ball{n}$ by $kn^{1/3}$ 
in the direction of the positive $x$-axis. This shifted set will be denoted by $\mathcal{B}_n^t$, which
is the set of all $(x,y)\in\C$ with $0\leq x\leq 2kn^{1/3}$ and
\begin{align*}
\abs{y} & \leq \frac{x^2}{3}, \quad \text{ for } 0\leq x\leq kn^{1/3}, \\
\abs{y} & \leq \frac{(2kn^{1/3}-x)^2}{3}, \quad  \text{ for }kn^{1/3}<x\leq 2kn^{1/3}.
\end{align*}
On the shifted set we define the function $g_n^t:\ball{n}^t\to\R$, by
\begin{equation}\label{eq:initial_translated_g}
g^t_n(x,y)=g_n(x + kn^{1/3}, y),
\end{equation}
which solves the same Dirichlet problem \eqref{eq:g_laplace} on $\ball{n}^t$
with the origin moved to $(kn^{1/3},0)$.
By symmetry of $g_n$, it is enough to compute $g_n^t$ for vertices $(x,y)$
with $0\leq x \leq k n^{1/3}$ and $y\geq 0$.
For $z=(x,y)\in\ball{n}^t$, with $y\neq 0$, the Laplace $\Delta g_n^t(z)$ is
equal to $-1/2$, hence on each ``tooth'' of the comb, $g^t_n$ satisfies the
linear recursion
\begin{equation*}
2 g^t_n(x,y) = g^t_n(x,y+1) + g^t_n(x,y-1) + 1,
\end{equation*}
which has the general solution
\begin{equation}
\label{eq:g_n_solution}
g^t_n(x,y) = \frac{1}{2}(y-y^2) + c_1(x) + y c_2(x),
\end{equation}
where $c_1(x)$ and $c_2(x)$ are functions of $x$, to be determined. For $(x,0),(x,1)\in\C$, we have
\begin{equation}
\label{eq:c1_c2}
g^t_n(x,0)=c_1(x) \quad \text{and}\quad g^t_n(x,1)=c_1(x)+c_2(x).
\end{equation}
From \eqref{eq:g_laplace} we have the boundary conditions $g_n^t(0,0)=0$ and $g_n^t(2kn^{1/3},0)=0$
and for $0\leq x\leq kn^{1/3}$, we have $ g_n^t(x,x^2/3)=0$.
On the other hand, from equation \eqref{eq:g_n_solution}, we get
\begin{equation*}
g_n^t(x,x^2/3)=\frac{x^2}{6}\Big(1-\frac{x^2}{3}\Big)+c_1(x)+\frac{x^2}{3}c_2(x)=0,
\end{equation*}
which implies that the function $c_2(x)$ can be written as
\begin{equation}
\label{eq:c_2}
c_2(x)=\frac{1}{2}\Big(\frac{x^2}{3}-1\Big)-\dfrac{3}{x^2}c_1(x). 
\end{equation}
Moreover, on the $x$-axis the Laplace operator of $g_n^t$ satisfies
\begin{align}\label{eq:laplace_x_axis}
\begin{split}
 \Delta g_n^t(x,0)& =
\begin{cases}
-\frac{1}{4}, & \text{ if } x\neq kn^{1/3} \\
-\frac{1}{4} (n+1) , & \text{ if } x=kn^{1/3}.
\end{cases}
\end{split}
\end{align}
For $x\neq kn^{1/3}$, that is, when $(x,0)$ is not the center of $\ball{n}^t$, we have
\begin{equation*}
 g_n^t(x+1,0)=4g_n^t(x,0)-g_n^t(x-1,0)-2g_n^t(x,1)-1,
\end{equation*}
and using \eqref{eq:c1_c2} we obtain
\begin{equation*}
 c_1(x+1)=2c_1(x)-c_1(x-1)-2c_2(x)-1,
\end{equation*}
which together with \eqref{eq:c_2} gives an equation for $c_1$, namely
\begin{equation*}
 c_1(x+1)=\Big(2+\frac{6}{x^2}\Big)c_1(x)-c_1(x-1)-\frac{x^2}{3}.
\end{equation*}
This has an explicit solution as a polynomial of degree $4$, given by
\begin{equation}
\label{eq:c_1}
 c_1(x)=-\frac{1}{18}x^4+bx^3-\frac{1}{36}x^2,
\end{equation}
where $b$ is a free parameter which can be computed using the other boundary conditions
for $g_n^t$. Since $\Delta g_n^t(kn^{1/3},0)=-\frac{1}{4}(n+1)$, using equations
\eqref{eq:c1_c2}, \eqref{eq:c_2}, and \eqref{eq:c_1}, we obtain
\begin{equation*}
b= \dfrac{5K+27n}{18(1+3K^2)},
\end{equation*}
where $K=kn^{1/3}$, and the constant $k= \left(\frac{3}{2}\right)^{2/3}$ is the
same as in Theorem \ref{thm:idla_cluster}. Since we are interested in the form of $g_n^t$
for $n$ sufficiently large, we expand $b$ around $n = \infty$, giving
\begin{equation*}
 b(n)=\frac{1}{6l}n^{1/3}+\O(n^{-1/3}).
\end{equation*}
Putting everything together we get $g_n(x,y) = g_n^t\big(k n^{1/3} - \abs{x}, \abs{y}\big)$, with
\begin{equation*}
g_n^t(x,y)=\Big(\frac{1}{6l}n^{1/3}+\O(n^{-1/3})\Big)(x^3-3xy) +\frac{1}{36}(3y-18y^2-2x^4-x^2 + 12 x^2 y).
\end{equation*}

\subsection{IDLA inner bound}
We are now able to conclude the proof of Theorem \ref{thm:idla_cluster}.
\begin{lem}
\label{lem:lambda_choice}
For all $\varepsilon > 0$ there exists $n_\varepsilon$, such that for all $n\geq n_\varepsilon$
and all $z\in \ball{n(1-\varepsilon)}$
\begin{equation*}
\frac{\varepsilon}{4}\leq \frac{\mathbb{E}[M - \tilde{L}]}{\mathbb{E}[M + \tilde{L}]}.
\end{equation*}
\end{lem}
\begin{proof}
By \eqref{eq:lambda}, one needs to study the function $\lambda_n(x,y) = \frac{f_n(x,y)}{g_n(x,y)}$.
We have
\begin{equation*}
\lambda_n(x,y) =\frac{\left(|y| - \frac{n_x}{2}\right)^2}{2 c_1(x) + \big(2c_2(x)+1\big) y - y^2},
\end{equation*}
where $c_1(x)$, $c_2(x)$ and $n_x$ are defined in \eqref{eq:c_1}, \eqref{eq:c_2} and \eqref{eq:n_x}, respectively.
It suffices to consider the first quadrant.
For every fixed $x$, the function $\lambda_n(x,y)$ is decreasing in $y$ for
$0 \leq y \leq \frac{n_x}{2}$. From the proof of Theorem \ref{thm:sandpile_cluster} we already know that
\begin{equation}
\label{eq:nx_expansion}
\frac{n_x}{2} = l\left(n^{1/3}-\frac{x}{k}\right)^2 + \O(1).
\end{equation}
For $0<\varepsilon<1$ consider the set
\begin{equation*}
\ball{n,\varepsilon} = \bigg\{(x,y)\in\C:\:  \abs{x}\leq (1-\varepsilon)k n^{1/3} \text{ and } 
                             \abs{y}\leq (1-\varepsilon)l \left(n^{1/3}-\frac{\abs{x}}{k}\right)^2\bigg\}.
\end{equation*}
Obviously $\ball{n,\varepsilon} \subset \ball{n}$ for all $\varepsilon$, hence $\frac{f_n}{g_n}$
is well defined on this set. Furthermore, by \eqref{eq:nx_expansion}, $\frac{f_n}{g_n}$ is also decreasing
on $\ball{n,\varepsilon}$ as a function of $y$, for all $\varepsilon>0$ and $n$ big enough. This
means that it is enough to study $\frac{f_n}{g_n}$ at the inner boundary of $\ball{n,\varepsilon}$.
For this, let $z=(x,y)\in\ball{n,\varepsilon}$ with $\abs{y} = (1-\varepsilon)l \left(n^{1/3}-\frac{\abs{x}}{k}\right)^2$ 
be such a boundary point. Then 
\begin{equation*}
\lim_{n\to\infty} \frac{f_n(z)}{g_n(z)} = \frac{\varepsilon}{4-\varepsilon} > \frac{\varepsilon}{4}.
\end{equation*}
The statement follows from the fact that for each  $\varepsilon>0$ one can find an $\varepsilon'>0$
such that $\ball{n(1-\varepsilon)} \subset \ball{n,\varepsilon'}$.
\end{proof}
\begin{proof}[Proof of Theorem \ref{thm:idla_cluster}]
Recall that we need to show the convergence of the series
\eqref{eq:borel_cantelli_idla}. Fix $z\in\ball{n(1-\varepsilon)}$. We set $\lambda=\frac{\varepsilon}{4}>0$ in
Lemma \ref{lem:lambda_choice}, and choose
\begin{equation*}
a=(1+\lambda)\mathbb{E}[\tilde{L}]=\Big(1+\frac{\varepsilon}{4}\Big)\mathbb{E}[\tilde{L}]
\end{equation*}
in equation \eqref{eq:ineq_a}.  Apply now Lemma \ref{lem:large_dev} to $M$ and $\tilde{L}$.
Recall also that $\mathbb{E}[M]>\mathbb{E}[\tilde{L}]$. Then
\begin{align*}
\begin{split}
\mathbb{P}[M\leq a]+\mathbb{P}[\tilde{L}\geq a] & =
\mathbb{P}\Big[M\leq \Big(1+\frac{\varepsilon}{4}\Big)
\mathbb{E}[\tilde{L}]\Big]+\mathbb{P}\Big[\tilde{L}
\geq \Big(1+\frac{\varepsilon}{4}\Big)\mathbb{E}[\tilde{L}]\Big]\\
& \leq 4 \exp\bigl\{-c_{\lambda}\mathbb{E}[\tilde{L}]\bigr\} \leq 4\exp\Bigl\{-c_{\lambda}\frac{g_n(z)-f_n(z)}{G_n(z,z)}\Bigr\},
\end{split}
\end{align*}
where $c_\lambda$ is a constant depending only on $\lambda$. Hence, for all $n\geq n_{\varepsilon}$, we have
\begin{equation}\label{eq:idla_comb_series}
\sum_{n\geq n_{\varepsilon}} \sum_{z\in \ball{n(1-\varepsilon)}}\mathbb{P}[z\notin A_n]
\leq 4 \sum_{n\geq n_{\varepsilon}} \sum_{z\in \ball{n(1-\varepsilon)}} \exp\Bigl\{-c_{\lambda}\frac{g_n(z)-f_n(z)}{G_n(z,z)}\Bigr\}.
\end{equation} 
In order to estimate the stopped Green function $G_n(z,z)$ upon exiting $\ball{n}$, with $z=(x,y)$, note that
\begin{equation*}
|y|\leq b_n(x):=l\Big(n^{1/3}-\frac{|x|}{k}\Big)^2.
\end{equation*}
We have the trivial upper bound $G_n(z,z) \leq 2 G_A(y,y)$ where $G_A$ is the Green function
of the simple random walk on the integer line, stopped at the interval $A=\big[-b_{n}(x),b_{n}(x)\big]$.
Using Proposition 1.6.3 and Theorem 1.6.4 from
\textsc{Lawler} \cite{lawler_intersections}, this can be bounded by
\begin{equation}\label{eq:g_a_order}
G_A(y,y)= \dfrac{b_{n}(x)^2-y^2}{b_{n}(x)}
\leq l\Big(n^{1/3}-\frac{|x|}{k}\Big)^2.
\end{equation}
For every $\varepsilon>0$, the function $g_n(z)-f_n(z)$ is again decreasing on every
non-crossing path which starts at $o$ and stays inside $\ball{n(1-\varepsilon)}$.
Hence, it attains its minimum on the inner boundary $\partial_I \ball{n(1-\varepsilon)}$
 of $\ball{n(1-\varepsilon)}$.
Taking limits, we get for every sequence $z_n=(x,y_n)$ with $x$ fixed and
$z_n\in\partial_I \ball{n(1-\varepsilon)}$
\begin{equation*}
\lim_{n\to\infty}\frac{g_n(z_n)-f_n(z_n)}{n^{4/3}} = \frac{k}{4} (2 - \varepsilon) \varepsilon,
\end{equation*}
and for the sequence $z'_n = (x_n, 0)$ with $x_n = k n^{1/3}(1-\varepsilon)^{1/3}$
\begin{equation*}
\lim_{n\to\infty}\frac{g_n(z'_n)-f_n(z'_n)}{n^{4/3}} =
\frac{k}{4} \big(3-2\varepsilon -(\varepsilon -3)(1-\varepsilon)^{1/3}\big).
\end{equation*}
Hence for all $\varepsilon>0$ and $n$ big enough
\begin{equation*}
\min_{z\in \ball{n(1-\varepsilon)}} \big(g_n(z)-f_n(z)\big) \geq C_{\varepsilon} \cdot n^{4/3},
\end{equation*}
for a constant $C_\varepsilon$ which depends only on $\varepsilon$.
Since, by \eqref{eq:g_a_order} the stopped Green function $G_A(z,z)$ is of order $\O(n^{2/3})$,
this implies
\begin{equation*}
\min_{z\in \ball{n(1-\varepsilon)}} \frac{g_n(z)-f_n(z)}{G_n(z,z)}\geq C'_{\varepsilon} \cdot n^{2/3}.
\end{equation*}
Hence, \eqref{eq:idla_comb_series} can be bounded by
\begin{align*}
 \sum_{n\geq n_{\varepsilon}} \sum_{z\in \ball{n(1-\varepsilon)}}\mathbb{P}[z\notin A_n]
&\leq 4 \sum_{n\geq n_{\varepsilon}} n \exp\{-c_{\lambda}C'_{\varepsilon}n^{2/3}\}<\infty,
\end{align*}
which concludes the proof.
\end{proof}

\subsection{The recurrent potential kernel}
In the recurrent lattice case $\Z^2$, the limiting shape of IDLA is derived using estimates
for the {\em recurrent potential kernel} which is defined as follows
\begin{equation*}
A(x,o) = \lim_{n\to\infty} \sum_{t=0}^n\big(\P_o[X_t = o] - \P_x[X_t = o]\big), \quad \text{for } x\in\C.
\end{equation*}
See \textsc{Lawler, Bramson and Griffeath} \cite{lawler_bramson_griffeath} for more details. For some constant $N>0$,
the {\em level sets of the potential kernel} are sets of the form $\{ x \in \C:\: A(x,o) \geq N\}$, 
and the {\em level sets of the Green function} are of the form $\{ x \in \C:\: G(x,o) \geq N\}$. 

In all previously known cases, the limiting shape of the IDLA-cluster is determined by
the level sets of the potential kernel in the recurrent case and by
level sets of the Green function in the transient case. Nevertheless, this is not the case for comb lattices
$\C$, even if the simple random walk is recurrent. In order to show this, let us consider the generating 
function of the potential kernel
\begin{align*}
A(x,o|z) &= \sum_{t=0}^\infty\big(\P_o[X_t = o] - \P_x[X_t = o]\big)z^t \\
         &= G(o,o|z) - G(x,o|z),
\end{align*}
where $G(x,y|z) = \sum_{t=0}^\infty \P_x[X_t = y] z^t$ is the generating function
of the Green function of simple random walk on $\C$.
Using standard techniques for generating functions one gets for $x = (x_1,x_2)\in \C$
\begin{equation*}
G(x,o|z) = F_1(z)^{\abs{x_2}} F_2(z)^{\abs{x_1}} G(o,o|z)
\end{equation*}
with
\begin{align*}
F_1(z) = \frac{1-\sqrt{1-z^2}}{z}, \qquad
F_2(z) = \frac{1 + \sqrt{1-z^2} - \sqrt{2}\sqrt{1-z^2+\sqrt{1-z^2}}}{z}
\end{align*}
and
\begin{align*}
G(o,o|z) = \frac{\sqrt{2}}{\sqrt{1-z^2 + \sqrt{1-z^2}}}.
\end{align*}
See \textsc{Bertacchi and Zucca} \cite{bertacchi_zucca} for details.
Therefore
\begin{equation*}
A(x,o|z) = G(o,o|z)\left(1 - F_1(z)^{\abs{x_2}} F_2(z)^{\abs{x_1}}\right) 
\end{equation*}
and $\lim_{z\to1^-} A(x,o|z) = 2 \abs{x_1}$.
By Abel's theorem for power series, the potential kernel on the comb $\C$, if it exists, has to 
be equal to $A(x,o) = 2 \abs{x_1}$. Hence we have an example of an IDLA cluster whose behaviour
is not given by the level sets of the potential kernel.
\section{Rotor-Router Aggregation}\label{sec:rotor_router}

A \emph{rotor-router walk} on a graph $G$ is defined as follows. For each vertex $x$ fix
a cyclic ordering $c(x)$ of its neighbours, i.e., $c(x)=(x_0,x_1,\ldots,x_{d(x)-1})$,
where  $x\sim x_i$ for all $i=0,1,\ldots ,d(x)-1$. The ordering $c(x)$ is called the
\emph{rotor sequence} of $x$. A \emph{rotor configuration} is a function
$\rho: G\to G$, with $\rho(x) \sim x$, for all $x\in G$. Hence $\rho$ assigns to every
vertex one of its neighbours.  A \emph{particle configuration} is a function
$\sigma:G\to \N_0$, with finite support. If $\sigma(x)=m>0$, we say that there are
$m$ particles at vertex $x$. A particle located at a vertex $x$ with current rotor
$\rho(x) = x_i$,  performs a rotor-router walk like this: it first sets $\rho(x) = x_{i+1}$, where addition
is modulo $d(x)$, and then it moves to $x_{i+1}$.

\emph{Rotor-router aggregation} is a deterministic process of increasing subsets
$(R_i)_{i\in\N}$ of $G$ defined recursively as $R_1 = \lbrace o \rbrace$, and
\begin{equation*}
R_{i} = R_{i-1} \cup \lbrace z_i \rbrace \quad\text{ for } i \geq 2,
\end{equation*}
where $z_i$ is the first vertex outside of $R_{i-1}$ that is visited
by a particle performing a rotor-router walk, started at $o$. The particle stops at $z_i$,
and a new particle starts its tour at the origin, but without resetting the rotor configuration.
The set $R_n$ of the occupied sites in $G$ is called
the \emph{rotor-router cluster} of $n$ particles.
The \emph{odometer function} $u_R(x)$ of rotor-router aggregation
is defined as the total number of particles which are sent out
by the vertex $x$ during the creation of the rotor-router cluster $R_n$.

In this section we study rotor-router aggregation on $\C$, and we give an inner bound for the cluster $R_n$ which
holds for arbitrary initial configuration of rotors and is independent of the rotor sequence.
The approach below relies on an idea of \textsc{Holroyd and Propp} \cite{holroyd_propp},
who use rotor weights in order to prove a variety of inequalities concerning rotor-walks and random walks.

\subsection{Rotor Weights}
Let $G$ be a locally finite and connected graph, $\sigma_0:G\to\N$ an initial particle
configuration with finite support, and $\rho_0: G\to G$ an initial rotor
configuration with $\rho_0(x) = x_0$ for all $x\in G$, that is, all
initial rotors point to the first element in the rotor sequence $c(x)$.
Routing particles in the system, such that at each time step $t$ exactly one particle
makes one step of a rotor-router walk, gives rise to a
sequence $(\rho_t, \sigma_t)_{t\in\N_0}$ of rotor- and particle-configurations.
To each of the possible states $(\rho_t,\sigma_t)$ of the system, we will assign a weight.

Fix a function $h:G\to \R$. Define the \emph{particle weights} at time $t$ to be
\begin{equation}
\label{eq:def_particles_weights}
\mathbf{W_P}(t) = \sum_{x\in G}\sigma_t(x) h(x),
\end{equation}
and the \emph{rotor weights} of single vertices $x\in G$ as
\begin{equation}
\label{eq:def_rotor_weights}
w(x, k) =
\begin{cases}
0, &\text{for } k = 0 \\
w(x, k-1) + h(x) - h\big(x_{k\bmod d(x)}\big), &\text {for } k > 0,
\end{cases}
\end{equation}
where $x_i$ is the $i$-th neighbour of $x$ in the rotor sequence $c(x)$.
Notice that, for $k \geq d(x)$,
\begin{equation}
\label{eq:rotor_weight_laplace}
w(x,k) = w\big(x,k-d(x)\big) - d(x)\Delta h(x).
\end{equation}
The total \emph{rotor weights} at time $t$ are given by
\begin{equation*}
\mathbf{W_R}(t) = \sum_{x\in G} w(x,u_t(x)),
\end{equation*}
where $u_t(z)$ is the odometer function of this process, that is, the number of
particles sent out by $x$ in the first $t$ steps. Note that $\rho_0$ is chosen in such a 
way that, if $i \equiv u_t(x) \bmod d(x)$, then $x_i = \rho_t(x)$ for all $t\geq 0$ and $x\in G$.

It is easy to check that the sum of particle- and rotor-weights are invariant under
routing of particles, i.e., for all times $t,t'\geq 0$
\begin{equation}
\label{eq:rotor_weights_invariance}
\mathbf{W_P}(t) + \mathbf{W_R}(t) = \mathbf{W_P}(t') + \mathbf{W_R}(t').
\end{equation}

For rotor-router aggregation on $G$ we start with $n$ particles at
the origin, that is, $\sigma_0 = n\cdot\delta_o$, and we route a particle
only if there is at least one other particle at the same position. The process terminates
when no two particles are at the same position. Denote by $t^\star = t^\star(n)$ the
number of steps it takes to finish the process, and by $(\sigma_{t^\star},\rho_{t^\star})$ the
final configuration. By the abelian property, the configuration $(\sigma_{t^\star},\rho_{t^\star})$ does not depend
on the order the particles made their steps, and by definition $\sigma_{t^\star}(x) = \indicator{\{x\in R_n\}}$.

We use the following weight function: for some $y\in G$, let $h_y:G\to\R$ given by
\begin{equation}
\label{eq:def_hy}
h(x) = h_y(x) = \frac{G_n(y,x)}{d(x)},
\end{equation}
where $G_n$ is the Green function of the simple random walk on $G$, stopped upon exiting the sandpile 
cluster $S_n$ with initial mass distribution $\mu_0=n\cdot\delta_{o}$ on $G$. Take now
some $y\in S_n$. Here $S_n$ is the sandpile cluster of the divisible sandpile
on some general graph $G$.
Recall that the Laplace of $h_y(x)$ on $G$ is given by \eqref{eq:laplace_hy}. The particle weights at the beginning are
\begin{equation}
\label{eq:pw_0}
\mathbf{W_P}(0) = n h_y(o),
\end{equation}
while the rotor weights are $\mathbf{W_R}(0) = 0$. At the end of the process, i.e., at time 
$t^\star$ when the rotor-router cluster $R_n$ is formed, we have
\begin{equation}
\label{eq:pw_end}
\mathbf{W_P}(t^\star) = \sum_{x\in R_n} h_y(x) \leq \sum_{x\in S_n} h_y(x),
\end{equation}
since $h_y$ is equal to $0$ outside of $S_n$. For the rotor weights we get from \eqref{eq:rotor_weight_laplace}
\begin{equation}
\label{eq:rw_end1}
\mathbf{W_R}(t^\star) = \sum_{x\in R_n} \left\lfloor \frac{u_R(x)}{d(x)} \right\rfloor
                             \big(-d(x)\Delta h_y(x)\big) + \sum_{x\in R_n} w(x, k_x),
\end{equation}
where $u_R$ is the rotor odometer function and $k_x = u_R(x) \bmod d(x)$.
By \eqref{eq:laplace_hy} and \eqref{eq:def_rotor_weights} 
\begin{equation}
\label{eq:rw_end}
\begin{split}
\mathbf{W_R}(t^\star) &= \left\lfloor \frac{u_R(y)}{d(y)} \right\rfloor +
   \sum_{x\in R_n} \sum_{i=0}^{k_x} \big(h_y(x) - h_y(x_i)\big) \\
   &\leq \frac{u_R(y)}{d(y)} + \sum_{x\in S_n}\sum_{z\sim x} \abs{h_y(x) - h_y(z)}.
\end{split}
\end{equation}
Hence by the invariance of the total weights \eqref{eq:rotor_weights_invariance}, we obtain
\begin{equation}
\label{eq:rotor_weight_inequality}
\sum_{x\in S_n} \big(n \delta_0(x) - 1\big) h_y(x) \leq \frac{u_R(y)}{d(y)} + \sum_{x\in S_n}\sum_{z\sim x} \abs{h_y(x) - h_y(z)}.
\end{equation}
Denote by $v(y)$ the lefthand side of \eqref{eq:rotor_weight_inequality}. Then $v(y)$
solves the Dirichlet problem
\begin{equation*}
\begin{cases}
\Delta v(y) = \frac{1}{d(y)}\big(1 - n \delta_0(y)\big), \text{ for } y \in S_n, \\
v(y) = 0, \text{ for } y \not\in S_n.
\end{cases}
\end{equation*}
By \eqref{eq:sandpile_odometer_laplace}, the normalized odometer function $u_n$ of
the divisible sandpile on $G$ with initial mass distribution $\mu_0 = n\cdot\delta_o$
satisfies exactly the same Dirichlet problem, hence $v(y) = u_n(y)$ and we get
the following result, which compares the odometer function $u_R$ of rotor-router
aggregation with the odometer function $u_n$ of the divisible sandpile. This result
holds for any locally finite and connected graph $G$.
\begin{prop}
\label{prop:sandpile_odo_leq_rotor_odo}
Let $u_n$ be the normalized odometer function of the divisible sandpile
with initial mass distribution $\mu_0 = n\cdot\delta_o$, and $u_R$
the odometer function of rotor-router aggregation with $n$ particles
starting at the origin $o\in G$. Then, for all $y\in G$,
\begin{equation}
\label{eq:sandpile_odo_leq_rotor_odo}
u_n(y) \leq \frac{u_R(y)}{d(y)} + \mathbf{\widetilde{W}}_\mathbf{R}(y),
\end{equation}
with
\begin{equation}
\label{eq:w_rest}
\mathbf{\widetilde{W}}_\mathbf{R}(y) =
  \sum_{x\in S_n}\sum_{z\sim x} \left|\frac{G_n(y,x)}{d(x)} - \frac{G_n(y,z)}{d(z)}\right|.
\end{equation}
\end{prop}
 \textsc{Levine and Peres}
\cite{peres_levine_strong_spherical} derived an inequality similar to \eqref{eq:sandpile_odo_leq_rotor_odo}
in the case of $\Z^d$ using a different method. For trees, $\mathbf{\widetilde{W}}_\mathbf{R}(y)$
can be expressed in terms of the expected distance from the starting point of a
random walk to the point where it first exits $S_n$.
\begin{prop}
\label{prop:wrest}
If $G$ is a tree and $\d(\cdot,\cdot)$ is the graph distance on $G$, then
\begin{equation*}
\mathbf{\widetilde{W}}_\mathbf{R}(x) = 2 \E_x\big[\d(x,X_T)\big] - 2,
\end{equation*}
where $T = \inf\big\{t\geq 0:\: X_t \not\in S_n\big\}$, and $(X_t)$ is the simple random walk on $G$.
\end{prop}
\begin{proof}
For $y\sim z$ let $N_{yz}$ be the number of transitions from $y$ to $z$
before the random walk exits $S_n$. Then
\begin{equation*}
\E_x\big[N_{yz} - N_{zy}\big] = \frac{G_n(x,y)}{d(y)} - \frac{G_n(x,z)}{d(z)}.
\end{equation*}
See also \cite[Proposition 2.2]{LP:book} for more details.
Since $G$ is a tree, the net number of crossings of each edge is smaller or equal to one, i.e.,
\begin{equation*}
\left|\E_x\big[N_{yz} - N_{zy}\big]\right| \leq 1.
\end{equation*}
We consider $G$ as a tree rooted at $x$, and denote by $\pi_{x,z}$ the shortest path from $x$ to $z$.
For $y\not= x$, write $y^-$ for the parent of $y$, i.e., the unique neighbour of $y$ that lies on
the shortest path $\pi_{x,y}$.
With this notation we get
\begin{equation*}
\sum_{\substack{y,z\in S_n\\y\sim z}} \left| \frac{G_n(x,y)}{d(y)} - \frac{G_n(x,z)}{d(z)} \right| =
\sum_{\substack{y,z\in S_n\\y\sim z}} \left| \E_x\big[N_{yz}-N_{zy}\big] \right| 
= 2 \sum_{\substack{y\in S_n\\y\not= x}} \E_x\big[N_{y^-y} - N_{yy^-}\big],
\end{equation*}
where the last equality is due to the antisymmetry of $N_{yz}-N_{zy}$. Let
\begin{equation*}
C_y = \big\{z\in S_n:\: y\in\pi_{x,z}\big\}
\end{equation*}
be the cone of $y$. The random variable $N_{y^-y} - N_{yy^-}$ is either zero or one, the latter if
the random walk exits $S_n$ in the cone $C_y$, hence
\begin{equation*}
\mathbf{\widetilde{W}}_\mathbf{R}(x) = 2 \sum_{\substack{y\in S_n\\ y\not= x}} \P_x\big[X_T\in C_y\big]
= 2 \sum_{\substack{y\in S_n\\ y\not= x}}\sum_{z\in C_y} \P_x\big[X_T=z\big].
\end{equation*}
For all $z\in\partial S_n$ we have
$\#\big\{y\in S_n\setminus\{x\}:\:z\in C_y\big\} = \d(x,z) - 1$, therefore
\begin{equation*}
\mathbf{\widetilde{W}}_\mathbf{R}(x) = 2 \sum_{z\in\partial S_n} \P_x\big[X_T = z\big]\big(\d(x,z) - 1\big) 
= 2 \E_x\big[\d(x,X_T)\big] - 2,
\end{equation*}
which completes the proof.
\end{proof}

\subsection{Rotor-router Aggregation on the Comb}

Since $\C$ is a tree, by Proposition \ref{prop:sandpile_odo_leq_rotor_odo} and \ref{prop:wrest}, 
one needs an upper bound for the expected distance  
from the starting point of a random walk $(X_t)$ on $\C$ to the point where it first exits $S_n$,
in order to derive an inner estimate of the rotor-router cluster.
Recall that on $\C$, the sandpile cluster $S_n$ has the shape $\ball{n}$ given in Theorem
\ref{thm:sandpile_cluster}. Using the trivial upper estimate
\begin{equation}
\label{eq:ez_upperbound}
\E_z\big[\d(z,X_T)\big]  \leq \max\big\{\d(z,w): w\in\partial S_n\big\} = \abs{x}+\abs{y} + l n^{2/3},
\end{equation}
with $z = (x,y)$ and $l = \frac{1}{2}\big(\frac{3}{2}\big)^{1/3}$ as in
Theorem \ref{thm:sandpile_cluster}, we can show the following inner bound.
\begin{thm}
\label{thm:weak_rotor_inner_bound}
Let $R_n$ be the rotor-router cluster of $n$ particles on  $\C$. Then, for $n\geq n_0$ and
for any initial rotor configuration and choice of rotor sequence, we have 
\begin{equation*}
\tilde{\mathcal{B}}_n \subset R_n,
\end{equation*}
where
\begin{equation*}
\tilde{\mathcal{B}}_n = \Big\{(x,y)\in\C:\: \abs{x} \leq k n^{1/3} - c_1 n^{1/6}, 
                                        \abs{y} \leq l\left(n^{1/3}-\frac{x}{k}\right)^2 + c_2 x - c_3 n^{1/3}\Big\},
\end{equation*}
where
\begin{equation*}
k = \left(\frac{3}{2}\right)^{2/3},
\quad
l = \frac{1}{2}\left(\frac{3}{2}\right)^{1/3}
\end{equation*}
and $c_1$, $c_2$ and $c_3$ are constants.
\end{thm}
\begin{proof}
By the definition  of rotor-router aggregation, $\big\lbrace z\in\C: u_R(z) > 0 \big\rbrace \subset R_n$,
and by Proposition \ref{prop:sandpile_odo_leq_rotor_odo} together with Proposition
\ref{prop:wrest}, we have for vertices $z=(x,y)$
\begin{equation*}
\frac{u_R(z)}{d(z)} \geq u_n(z) - 2 \E_z\big[\d(z, X_T)\big] + 2
                    \geq u_n(z) - 2 \big(\abs{x}+\abs{y} + l n^{2/3}\big) + 2,
\end{equation*}
where the last inequality is due to \eqref{eq:ez_upperbound}. By Corollary
\ref{cor:sandpile_odometer}, we have a lower bound of the sandpile odometer
$u_n$ for $z\in S_n$
\begin{equation*}
\gamma_n(z) - a \leq u_n(z),
\end{equation*}
where $a$ is a positive constant smaller than $2$, and $\gamma_n$ is the
function defined in \eqref{eq:odometer_C2}. Thus, to derive an inner bound, 
it suffices to check for which $z=(x,y)\in S_n$ the inequality
\begin{equation}
\label{eq:rr_inner_bound_ineq}
\gamma_n(x,y) - 2 \big(\abs{x}+\abs{y} + l n^{2/3}\big) > 0
\end{equation}
holds. By symmetry it is enough to consider $x,y\geq 0$.
We first check inequality \eqref{eq:rr_inner_bound_ineq} on a ``tooth''
of the comb, that is, for a fixed $x$. The function $\gamma_n$ is given as
\begin{equation*}
\gamma_n(x,y) = \frac{1}{2}\left(y-\frac{n_x}{2}\right)^2,
\end{equation*}
where $n_x$ is the amount of mass that ends up in the $x$-``tooth'' of
the sandpile. Since $x$ is fixed, we can treat $n_x$ as a constant. Hence
the right hand side of \eqref{eq:rr_inner_bound_ineq} is a quadratic polynomial
in $y$ with smallest root
\begin{equation*}
y_x = 2 + \frac{n_x}{2} - \sqrt{4 + \frac{k}{l} n^{2/3} + 2 n_x + 4 x}.
\end{equation*}
Substituting $n_x$ as calculated in \eqref{eq:n_x}, an expansion around $n=\infty$ gives
\begin{equation*}
y_x = l n^{2/3} - \frac{1}{2l} n^{1/3} x + \frac{x^2}{3} + \frac{2+\sqrt{6}}{3}x + c_1 n^{1/3} - c_2\frac{x^4}{n}.
\end{equation*}
Since $(x,y)\in S_n$, we have the bound $x\leq k n^{1/3}$, hence
\begin{equation}
\label{eq:rotor_yx_bound}
y_x = l \left(n^{1/3} - \frac{x}{k}\right)^2 + \frac{2+\sqrt{6}}{3}x - c n^{1/3},
\end{equation}
for $n\geq n_0$, and a positive constant $c$.
To get a bound on the $x$-axis, we calculate for which $x>0$ the inequality $y_x > 0$ is satisfied.
Since $y_x$ is a polynomial of degree $2$ in $x$ this is easy to do, and again by series expansion
around $n=\infty$ we obtain
\begin{equation}
\label{eq:rotor_x_bound}
x \leq k\cdot n^{1/3} - c_3 n^{1/6},
\end{equation}
for $n\geq n_0$.
The inner bound for $R_n$ now follows from \eqref{eq:rotor_x_bound} together with \eqref{eq:rotor_yx_bound}.
\end{proof}

\begin{minipage}{0.65\linewidth}
Figure \ref{fig:rotor_inner_bound} shows the inner estimate of the rotor-router cluster 
from Theorem \ref{thm:weak_rotor_inner_bound} in comparison to sandpile cluster $S_n$, for $n=1000$. 
The white area is the area where the inequality \eqref{eq:rr_inner_bound_ineq} is valid, and corresponds to the set
$\tilde{\mathcal{B}}_n$ of Theorem \ref{thm:weak_rotor_inner_bound}.
The colouring is based on the value of the right-hand side of \eqref{eq:rr_inner_bound_ineq}.\\

The inner bound could be improved if one has a substantially better upper bound for $\E_z[\d(z,X_T)]$. 
For regular graphs, one can also give an universal inner estimate for rotor-router aggregation, 
which relates the rotor-router cluster to a divisible sandpile cluster with a smaller mass. 
Using the methods of \textsc{Levine and Peres} in \cite{peres_levine_strong_spherical} one can deduce the following.
 \begin{prop}
Let $G$ be a regular graph with degree $d$ and root $o$, and let $R_n$ be the
rotor-router cluster of $n$ particles starting at $o$. Further, let $S_n$ be
the divisible sandpile cluster with mass distribution $\mu_0(x) = n\cdot\delta_o(x)$. Then
$S_{n/(2d-1)} \subset R_n.$
\end{prop}
\end{minipage}
\hspace*{0.5cm}
\begin{minipage}{0.3\linewidth}
\centering
\includegraphics[height=9.5cm]{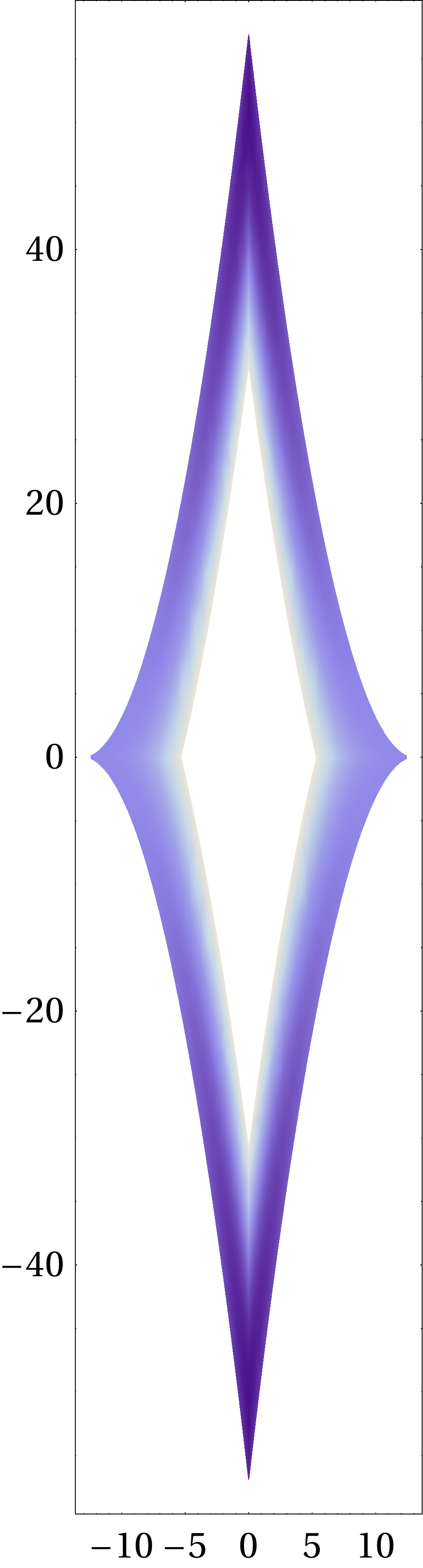}
\captionof{figure}{Inner bound}
\label{fig:rotor_inner_bound}
\end{minipage}

\paragraph*{Remarks.}
In the IDLA model, for an upper bound of the type $A_n\subset \ball{n(1+\varepsilon)}$, we still do not 
have sufficient tools. In all previous studied cases, the IDLA cluster grows uniformly, and this makes easy the study of random 
walks. This is of course violated in our case, since the set $\ball{n}$ in \eqref{eq:limit_shape}
grows with rate $n^{1/3}$ in the $x$-direction, and with rate $n^{2/3}$ in the $y$-direction. 
The \emph{harmonic measure} for random walks stopped when exiting the set $\ball{n}$, was studied in 
\cite{huss_sava_rotor}, by making use of a special rotor-router type process.
This is the first example where IDLA aggregate is not a set of uniform harmonic measure.
In \cite{huss_sava_rotor}, we also give subsets of $\C$ with uniform harmonic measure.

Recently \textsc{Duminil-Copin et al.} \cite{duminil_copin_idla} developed a method
for proving an outer bound for IDLA without needing a harmonic measure estimate, provided an inner bound is known.
Nevertheless, their method cannot be applied in our setting, since one of the required 
regularity conditions (weaker lower bound) is not satisfied on the comb.

\paragraph*{Acknowledgements.} The research of Wilfried Huss was partially 
supported by the FWF program P19115-N18 and the research of Ecaterina Sava
was supported by the FWF program W 1230-N13. Part of the work was done
during the stay of the first author at the University of Siegen, Germany.


\bibliographystyle{mypaperhep}

\vspace{1cm}
\begin{minipage}{0.45\textwidth}
Wilfried Huss\\
Vienna University of Technology\\
E-mail: \texttt{whuss@mail.tuwien.ac.at}\\
\texttt{http://www.math.tugraz.at/$\sim$huss}
\end{minipage}
\hfill
\begin{minipage}{0.45\textwidth}
Ecaterina Sava\\
Graz University of Technology\\
E-mail: \texttt{sava@tugraz.at}\\
\texttt{http://www.math.tugraz.at/$\sim$sava}
\end{minipage}

\end{document}

%% file: comb.tex
\begin{tikzpicture}
\foreach \x in {-2,...,2}
{
   \draw (\x, 2.5) -- (\x, -2.5);
   \foreach \y in {-2,...,2}
   \fill (\x,\y) circle (1.5pt);
}

\draw (-2.5, 0) -- (2.5, 0);

\coordinate[label=-45:$o$] (O) at (0,0);
\coordinate[label=0:{$z=(x,y)$}] (Z) at (2, 1);
\end{tikzpicture}

%% file: idla_comb.tex
\begin{tikzpicture}[scale=0.09]

\begin{scope}[xshift=-20cm]
\draw (-4,0) -- (-4,1);
\draw (-3,-3) -- (-3,3);
\draw (-2,-5) -- (-2,4);
\draw (-1,-9) -- (-1,8);
\draw (0,-11) -- (0,12);
\draw (1,-9) -- (1,10);
\draw (2,-4) -- (2,5);
\draw (3,-1) -- (3,2);
\draw (4,-1) -- (4,1);
\draw (5,0) -- (5,1);
\draw (6,0) -- (6,0);
\draw (-4,0) -- (6,0);
\end{scope}

\begin{scope}
\draw (-10,0) -- (-10,0);
\draw (-9,-2) -- (-9,1);
\draw (-8,-3) -- (-8,3);
\draw (-7,-2) -- (-7,5);
\draw (-6,-5) -- (-6,8);
\draw (-5,-11) -- (-5,8);
\draw (-4,-14) -- (-4,14);
\draw (-3,-17) -- (-3,19);
\draw (-2,-24) -- (-2,24);
\draw (-1,-26) -- (-1,26);
\draw (0,-40) -- (0,33);
\draw (1,-26) -- (1,29);
\draw (2,-19) -- (2,21);
\draw (3,-18) -- (3,18);
\draw (4,-12) -- (4,16);
\draw (5,-10) -- (5,8);
\draw (6,-4) -- (6,5);
\draw (7,-1) -- (7,2);
\draw (8,-1) -- (8,3);
\draw (9,-1) -- (9,2);
\draw (-10,0) -- (9,0);
\end{scope}

\begin{scope}[xshift=26cm]
\draw (-12,0) -- (-12,0);
\draw (-11,-1) -- (-11,1);
\draw (-10,-2) -- (-10,5);
\draw (-9,-5) -- (-9,5);
\draw (-8,-8) -- (-8,9);
\draw (-7,-10) -- (-7,12);
\draw (-6,-13) -- (-6,15);
\draw (-5,-23) -- (-5,24);
\draw (-4,-24) -- (-4,30);
\draw (-3,-31) -- (-3,33);
\draw (-2,-42) -- (-2,41);
\draw (-1,-45) -- (-1,47);
\draw (0,-60) -- (0,51);
\draw (1,-52) -- (1,48);
\draw (2,-35) -- (2,45);
\draw (3,-38) -- (3,33);
\draw (4,-25) -- (4,28);
\draw (5,-27) -- (5,22);
\draw (6,-13) -- (6,15);
\draw (7,-9) -- (7,12);
\draw (8,-9) -- (8,7);
\draw (9,-7) -- (9,5);
\draw (10,-3) -- (10,4);
\draw (11,0) -- (11,2);
\draw (12,0) -- (12,0);
\draw (-12,0) -- (12,0);
\end{scope}

\end{tikzpicture}